\documentclass[12pt,a4paper]{article}

\newtheorem{theorem}{Theorem}[section]
\newtheorem{proposition}[theorem]{Proposition}
\newtheorem{lemma}[theorem]{Lemma}
\newtheorem{corollary}[theorem]{Corollary}

\usepackage{hyperref}
\usepackage{amsmath}
\usepackage{amssymb}
\usepackage{array}
\usepackage[english]{babel}
\usepackage{calc}
\usepackage{enumerate}
\usepackage[latin1]{inputenc}
\usepackage{latexsym}

\begin{document}
\title{Compactness and weak-star continuity of derivations on weighted 
convolution algebras}
\author{Thomas Vils Pedersen}

\maketitle


\begin{abstract}
\noindent
Let $\omega$ be a continuous weight on $\mathbb R^+$ and let $L^1(\omega)$ be the corresponding 
convolution algebra. By results of Grønbæk and Bade \& Dales the continuous derivations from $L^1(\omega)$ to its dual space $L^{\infty}(1/\omega)$ are exactly the maps of the form
$$(D_{\varphi}f)(t)=\int_0^{\infty}f(s)\,\frac{s}{t+s}\,\varphi(t+s)\,ds
\qquad\text{($t\in\mathbb R^+$ and $f\in L^1(\omega)$)}$$
for some $\varphi\in L^{\infty}(1/\omega)$. Also, every $D_{\varphi}$ has a unique extension to a continuous derivation $\overline{D}_{\varphi}:M(\omega)\to L^{\infty}(1/\omega)$ from the corresponding measure algebra. We show that a certain condition on $\varphi$ implies that $\overline{D}_{\varphi}$ is weak-star continuous. The condition holds for instance if $\varphi\in L_0^{\infty}(1/\omega)$. We also provide examples of functions $\varphi$ for which $\overline{D}_{\varphi}$ is not weak-star continuous. Similarly, we show that $D_{\varphi}$ and $\overline{D}_{\varphi}$ are compact under certain conditions on $\varphi$. For instance this holds if $\varphi\in C_0(1/\omega)$ with $\varphi(0)=0$. Finally, we give various examples of functions $\varphi$ for which $D_{\varphi}$ and $\overline{D}_{\varphi}$ are not compact.
\end{abstract}

\section{Introduction}
\label{sec:intro}

Traditionally the study of derivations from a Banach algebra to its Banach modules has mainly focused on the existence of such derivations. In some recent papers by Choi and Heath the aim has instead been to characterise the derivations from a concrete Banach algebra to its dual space, and then to use this characterisation to study properties of the derivations: Every bounded  derivation from $l^1(\mathbb Z^+)$ to its dual space is of the form
$$D_{\psi}(\delta_0)=0\qquad\text{and}\qquad 
D_{\psi}(\delta_j)(\delta_k)=\frac{j}{j+k}\,\psi_{j+k}\ 
(j,k\in\mathbb Z^+,\ j\neq0)$$
for some $\psi\in l^{\infty}(\mathbb Z^+)$. It was shown in \cite{Hea} that $D_{\psi}$ is compact if and only if $\psi\in c_0(\mathbb Z^+)$. Moreover, 
the weakly compact derivations from $l^1(\mathbb Z^+)$ to its dual space are described in \cite{Ch-He:Trans} in terms of so-called translation-finite sets. Finally, the compact derivations from the disc algebra to its dual space are characterised in \cite{Ch-He}. In this paper we continue this line of thinking and consider properties of derivations from weighted convolution algebras $L^1(\omega)$ on $\mathbb R^+$ to their dual spaces. 
 
Let $L^1(\mathbb R^+)$ be the Banach space of (equivalence classes of)
integrable functions $f$ on $\mathbb R^+=[0,\infty)$ with the norm 
$\|f\|=\int_0^{\infty}|f(t)|\,dt$.
Throughout this paper $\omega$ will be a continuous weight function on $\mathbb R^+$, that is, a positive and continuous function on $\mathbb R^+$ satisfying $\omega(0)=1$ and $\omega(t+s)\le\omega(t)\omega(s)$ for all $t,s\in\mathbb R^+$.
We then define $L^1({\omega})$ as the weighted space of functions $f$ on 
$\mathbb R^+$ 
for which $f\omega\in L^1(\mathbb R^+)$ with the inherited norm 
$$\|f\|=\int_0^{\infty}|f(t)|\omega(t)\,dt.$$
With the usual convolution product
$$(f\ast g)(t)=\int_0^t f(s)g(t-s)\,ds\qquad
\text{for $t\in\mathbb R^+$ and $f,g\in L^1({\omega})$},$$
it is well known that $L^1({\omega})$ is a commutative Banach algebra.
Similarly, the space $M(\omega)$ of locally finite, complex Borel measures 
$\mu$ on $\mathbb R^+$ for which 
$$\|\mu\|=\int_0^{\infty}\omega(t)\,d|\mu|(t)<\infty$$
is a Banach algebra under convolution and contains $L^1({\omega})$ as a closed 
ideal. Also, $M(\omega)$ can be identified as the multiplier algebra of $L^1(\omega)$ 
and this induces a strong topology on $M(\omega)$ by identifying a measure with the corresponding convolution operator.

Moreover, we let $L^{\infty}(1/\omega)$ denote the Banach
space of measurable functions $\varphi$ on $\mathbb R^+$ for which $\varphi/\omega$ is essentially bounded with the norm
$\|\varphi\|=\text{ess\,sup}_{t\in\mathbb R^+}\,|\varphi(t)|/\omega(t)$. It is well known that the duality
$\langle f,\varphi\rangle=\int_0^{\infty}f(t)\varphi(t)\,dt$ for 
$f\in L^1(\omega)$ and $\varphi\in L^{\infty}(1/\omega)$ identifies 
$L^{\infty}(1/\omega)$
isometrically isomorphically with the dual space of $L^1(\omega)$. 

We denote by $C_b(1/\omega)$ the closed subspace of $L^{\infty}(1/\omega)$ of continuous functions in $L^{\infty}(1/\omega)$, and by $C_0(1/\omega)$ the closed subspace of $C_b(1/\omega)$ of functions $h\in C_b(1/\omega)$ for which $h/\omega$ vanishes at infinity. 
Then $M(\omega)$ is isometrically isomorphic to the dual space of $C_0(1/\omega)$ with the 
duality being defined by
$$\langle h,\mu\rangle
=\int_0^{\infty}h(t)\,d\mu(t)\qquad\text{for $h\in C_0(1/\omega)$ and $\mu\in M(\omega)$}.$$
We will need yet another closed subspace of $L^{\infty}(1/\omega)$. 
For $\varphi\in L^{\infty}(\mathbb R^+)$ we say that $\varphi(t)\to0$ as $t\to\infty$ if
$$\text{ess\,sup}_{t\ge T}\,|\varphi(t)|\to0\qquad\text{as }T\to\infty.$$
Similarly, we say that $\varphi(t)\to0$ as $t\to0$ if
$\text{ess\,sup}_{t\le T}\,|\varphi(t)|\to0$ as $T\to0$.
We then define $L_0^{\infty}(1/\omega)$ to be the closed subspace of $L^{\infty}(1/\omega)$ of those 
$\varphi\in L^{\infty}(1/\omega)$ for which $\varphi(t)/\omega(t)\to0$ as $t\to\infty$.

Recall that the dual space $L^{\infty}(1/\omega)=L^1(\omega)^*$ becomes a Banach $L^1(\omega)$-module via the action
$$\langle f,g\cdot\varphi\rangle=\langle f*g,\varphi\rangle\qquad
\text{for $f,g\in L^1(\omega)$ and $\varphi\in L^{\infty}(1/\omega)$}.$$
An easy calculation shows that the module action can be expressed as
$$(g\cdot\varphi)(t)=\int_0^{\infty}g(s)\varphi(t+s)\,ds\qquad
\text{for $t\in\mathbb R^+,\ g\in L^1(\omega)$ and $\varphi\in L^{\infty}(1/\omega)$}.$$
In particular it follows that $C_0(1/\omega)$ is a Banach $L^1(\omega)$-submodule of $L^{\infty}(1/\omega)$.
Also, if we consider $M(\omega)=C_0(1/\omega)^*$ as a dual Banach $L^1(\omega)$-module, then 
\begin{align*}
\langle h,g\cdot\mu\rangle
&= \langle g\cdot h,\mu\rangle
=\int_0^{\infty}\int_0^{\infty}g(s)h(t+s)\,ds\,d\mu(t)\\
&= \int_0^{\infty}\int_t^{\infty}g(r-t)h(r)\,dr\,d\mu(t)\\
&= \int_0^{\infty}\int_0^rg(r-t)\,d\mu(t)\,h(r)\,dr
=\langle h,g*\mu\rangle
\end{align*}
for $g\in L^1(\omega),\ h\in C_0(1/\omega)$ and $\mu\in M(\omega)$.
The dual module action $g\cdot\mu$ of $g\in L^1(\omega)$ on $\mu\in M(\omega)$ thus coincides with the usual convolution product $g\ast\mu$.
(One may also have wished for the product $g\cdot\varphi$ for $g\in L^1(\omega)$ and $\varphi\in L^{\infty}(1/\omega)$ to coincide with the usual convolution product of $g$ and $\varphi$. This could be obtained by choosing instead to identify $L^1(\omega)^*$ with the space 
$L^{\infty}(\mathbb R^-,\omega(-t))$ on $\mathbb R^-$ with the duality 
$\langle f,\varphi\rangle=\int_0^{\infty}f(t)\varphi(-t)\,dt$ for $f\in L^1(\omega)$ and 
$\varphi\in L^{\infty}(\mathbb R^-,\omega(-t))$. This approach is taken in, for instance, \cite{Da:Book}. We prefer instead (as in \cite{Gro}) to represent all our spaces on $\mathbb R^+$, and pay the price of the form of the product 
$g\cdot\varphi$.)

We now turn to derivations from $L^1(\omega)$ to its dual space $L^{\infty}(1/\omega)$. Recall that a linear map 
$D:L^1(\omega)\to L^{\infty}(1/\omega)$ is called a derivation if
$$D(f*g)=f\cdot Dg+g\cdot Df\qquad\text{for }f,g\in L^1(\omega).$$
The main part of the following result was proved by Grønbæk (\cite[Theorem~3.7]{Gro}). Bade and Dales 
(\cite[Theorem~2.3]{Ba-Da:Cont} or \cite[Theorem~5.6.27]{Da:Book}) then elaborated on Grønbæk's result to obtain the following.
For $s\in\mathbb R^+$ we denote by $\delta_s$ the unit point measure at $s$.

\begin{theorem}[Grønbæk and Bade \& Dales] 
\label{th:gr}
Let $\varphi\in L^{\infty}(1/\omega)$. Then
$$(D_{\varphi}f)(t)=\int_0^{\infty}f(s)\,\frac{s}{t+s}\,\varphi(t+s)\,ds
\qquad\text{for $t\in\mathbb R^+$ and $f\in L^1(\omega)$}$$
defines a continuous derivation from $L^1(\omega)$ to $L^{\infty}(1/\omega)$. Moreover, $D_{\varphi}$ has a unique extension to a continuous derivation $\overline{D}_{\varphi}:M(\omega)\to L^{\infty}(1/\omega)$. Also, $\overline{D}_{\varphi}$ is continuous when $M(\omega)$ is equipped with its strong topology and $L^{\infty}(1/\omega)$ with its weak-star topology, and 
$$(\overline{D}_{\varphi}\delta_s)(t)=\frac{s}{t+s}\,\varphi(t+s)\qquad\text{for $t,s\in\mathbb R^+$}.$$
Conversely, every continuous derivation from $L^1(\omega)$ to $L^{\infty}(1/\omega)$ equals $D_{\varphi}$ for some $\varphi\in L^{\infty}(1/\omega)$.
\end{theorem}

We mention that if we let $X$ be the densely defined operator on $L^1(\omega)$ given by $(Xf)(t)=tf(t)$ for $t\in\mathbb R^+$ and suitable $f\in L^1(\omega)$, and similarly on $L^{\infty}(1/\omega)$, then we have $D_{\varphi}f=(Xf)\cdot(X^{-1}\varphi)$ for $f$ and $\varphi$ in dense subsets of $L^1(\omega)$ resp. $L^{\infty}(1/\omega)$.

In this paper we study various properties of the derivations $\overline{D}_{\varphi}$. In Section~\ref{sec:weak-star} we consider weak-star  continuity, and in Section~\ref{sec:various} we rely on some of the results from Section~\ref{sec:weak-star} to prove various range and continuity properties of the derivations. Finally, some of these results are used in Section~\ref{sec:compact}, where compactness of the derivations is investigated. We remark that most of the results in this paper also are of interest in the unweighted case where $\omega\equiv1$.

\section{Weak-star continuity}
\label{sec:weak-star}

In this section we will study the weak-star continuity of the derivations $\overline{D}_{\varphi}:M(\omega)\to L^{\infty}(1/\omega)$. This is inspired by a similar result for homomorphisms between weighted algebras due to Grabiner: Let $\omega_1$ and $\omega_2$ be weights, and let $\Phi:L^1(\omega_1)\to L^1(\omega_2)$ be a non-zero continuous homomorphism. Then $\Phi$ has a unique extension to a continuous homomorphism $\widetilde{\Phi}:M(\omega_1)\to M(\omega_2)$ (\cite[Theorems~3.4]{Gr}) and this extension is automatically weak-star  continuous (\cite[Theorem~1.1]{Gr-wks}). 
Also, similar results hold for homomorphisms from $L^1(\omega)$ into some other commutative Banach algebras (\cite{Pe:wkprop}). We also mention that it  easily can be seen that a bounded derivation $D_{\psi}$ from $l^1(\mathbb Z^+)$ to its dual space is weak-star continuous if and only if 
$\psi\in c_0(\mathbb Z^+)$. 

For $f,g\in L^1(\omega)$ and $\varphi\in L^{\infty}(1/\omega)$ it follows from Fubini's theorem that
\begin{eqnarray}
\label{eq:tfi}
\langle f,D_{\varphi}g\rangle
&=& \int_0^{\infty}f(s)\int_0^{\infty}\frac{t}{t+s}\,\varphi(t+s)g(t)\,dt\,ds
\nonumber\\
&=& \int_0^{\infty}\left(\int_0^{\infty}f(s)\,\frac{t}{t+s}\,\varphi(t+s)\,ds\right)
g(t)\,dt.
\end{eqnarray}
This leads us to the following definition. Let
$$(T_{\varphi}f)(t)=\int_0^{\infty}f(s)\,\frac{t}{t+s}\,\varphi(t+s)\,ds
=\langle f,\overline{D}_{\varphi}\delta_t\rangle$$
for $f\in L^1(\omega),\,\varphi\in L^{\infty}(1/\omega)$ and $t\in\mathbb R^+$.

\begin{proposition}
\label{pr:tfi}
Let $\varphi\in L^{\infty}(1/\omega)$. 
\begin{enumerate}[(a)]
\item 
$(\overline{D}_{\varphi}\delta_s)$ is weak-star continuous in $L^{\infty}(1/\omega)$ for $s\in\mathbb R^+$.
\item
$T_{\varphi}$ is a continuous linear operator 
$T_{\varphi}:L^1(\omega)\to C_b(1/\omega)$.
\end{enumerate}
\end{proposition}

\noindent{\bf Proof}\quad
(a):\quad
Translation is continuous in $L^1(\omega)$ (\cite[Lemma~4.7.6]{Da:Book}), so translation is weak-star  continuous in $L^{\infty}(1/\omega)$. Hence $(\overline{D}_{\varphi}\delta_s)$ is weak-star continuous in $L^{\infty}(1/\omega)$ for $s>0$. Also, for $f\in L^1(\omega)$ we have
$$|\langle f,\overline{D}_{\varphi}\delta_s\rangle|
=\left|\int_0^{\infty}f(t)\frac{s}{t+s}\,\varphi(t+s)\,dt\right|
\le\|\varphi\|\omega(s)\int_0^{\infty}|f(t)|\omega(t)\frac{s}{t+s}\,dt\to0$$
as $s\to0$ by Lebesgue's dominated convergence theorem, so 
$\overline{D}_{\varphi}\delta_s\to0=\overline{D}_{\varphi}\delta_0$ weak-star in $L^{\infty}(1/\omega)$ as $s\to0$. 

(b):\quad
Let $f\in L^1(\omega)$. The estimate
$$|(T_{\varphi}f)(t)|\le\|\varphi\|\int_0^{\infty}|f(s)|\omega(t+s)\,ds
\le\|\varphi\|\cdot\|f\|\omega(t)\qquad\text{for $t\in\mathbb R^+$}$$
shows that $T_{\varphi}$ defines a continuous linear operator $T_{\varphi}:L^1(\omega)\to L^{\infty}(1/\omega)$.
Moreover, it follows from (a) that $T_{\varphi}$ maps into $C_b(1/\omega)$.
{\nopagebreak\hfill\raggedleft$\Box$\bigskip}

We will need the next couple of results. For $n\in\mathbb N$ let $e_n=n\cdot 1_{[0,1/n]}$. It is well known that $(e_n)$ is a bounded approximate identity for $L^1(\omega)$. Also, $\langle h,\mu\rangle=\int_0^{\infty}h(t)\,d\mu(t)$ is well-defined for $h\in C_b(1/\omega)$ and $\mu\in M(\omega)$. 
 
\begin{lemma}
\label{le:bai}
\ 
\begin{enumerate}[(a)]
\item 
Let $h\in C_0(1/\omega)$. Then $e_n\cdot h\to h$ in $C_0(1/\omega)$ as $n\to\infty$.
\item 
Let $h\in C_b(1/\omega)$ and $\mu\in M(\omega)$. Then 
$\langle e_n\cdot h,\mu\rangle\to\langle h,\mu\rangle$ as $n\to\infty$.
\end{enumerate}
\end{lemma}

\noindent{\bf Proof}\quad
(a):\quad
For $n\in\mathbb N$ and $t\in\mathbb R^+$ we have
$$(e_n\cdot h-h)(t)=\int_0^{1/n}n(h(t+s)-h(t))ds.$$
Given $\varepsilon>0$, we choose $T\in\mathbb R^+$ such that $|h(t)|/\omega(t)<\varepsilon$ for $t\ge T$. For all $n\in\mathbb N$ we then have 
$|(e_n\cdot h-h)(t)|/\omega(t)<(1+\sup_{s\in[0,1]}\omega(s))\varepsilon$ for $t\ge T$. Since $h$ is uniformly continuous on $[0,T+1]$ we can choose $N\in\mathbb N$ such that 
$|h(t+s)-h(t)|<\varepsilon\cdot\inf_{t\le T}\omega(t)$ for all $0\le t\le T$ and $s\le\frac{1}{N}$. Hence 
$$\frac{|(e_n\cdot h-h)(t)|}{\omega(t)}
\le\sup_{s\le1/N}\frac{|h(t+s)-h(t)|}{\omega(t)}<\varepsilon$$
for all $0\le t\le T$ and $n\ge N$. This finishes the proof.

(b):\quad
Given $\varepsilon>0$, we choose $T\in\mathbb R^+$ such that $|\mu\cdot\omega|([T,\infty))<\varepsilon$. We have
\begin{align*}
|\langle e_n\cdot h-h,\mu|
&\le\int_0^T|(e_n\cdot h-h)(t)|\,d|\mu|(t)
+\int_T^{\infty}|(e_n\cdot h-h)(t)|\,d|\mu|(t)\\
&\le\|\mu\|\sup_{0\le t\le T}\frac{|(e_n\cdot h-h)(t)|}{\omega(t)}
+\|e_n\cdot h-h\|\cdot|\mu\cdot\omega|([T,\infty)).
\end{align*}
The second term is bounded by $C\varepsilon$ and it follows from the proof of part (a) that there exists $N\in\mathbb N$ such that the same holds for the first term for $n\ge N$.
{\nopagebreak\hfill\raggedleft$\Box$\bigskip}

\begin{corollary}
\label{co:bai}
Let $\mu\in M(\omega)$. Then $e_n\ast\mu\to\mu$ strongly in $M(\omega)$ and weak star in $C_b(1/\omega)^*$ (and in particular weak-star in $M(\omega)$) as $n\to\infty$.
\end{corollary}

\noindent{\bf Proof}\quad
For $f\in L^1(\omega)$ we have $e_n\ast\mu\ast f\to\mu\ast f$ in $L^1(\omega)$ as $n\to\infty$ (since $(e_n)$ is a bounded approximate identity for $L^1(\omega)$). Hence $e_n\ast\mu\to\mu$ strongly in $M(\omega)$ as $n\to\infty$. Moreover, for $h\in C_b(1/\omega)$ we have 
$\langle h,e_n\ast\mu\rangle=\langle e_n\cdot h,\mu\rangle
\to\langle h,\mu\rangle$ as $n\to\infty$ by Lemma~\ref{le:bai}(b). Hence $e_n\ast\mu\to\mu$ weak-star in $C_b(1/\omega)^*$ as $n\to\infty$.
{\nopagebreak\hfill\raggedleft$\Box$\bigskip}

The adjoint of a continuous linear operator is weak-star continuous. The following result shows that the converse also is true for the operators $\overline{D}_{\varphi}$.

\begin{proposition}
\label{pr:wkscts}
For $\varphi\in L^{\infty}(1/\omega)$ the following conditions are equivalent:
\begin{enumerate}[(a)]
\item 
$\overline{D}_{\varphi}$ is weak-star continuous.
\item
$\overline{D}_{\varphi}\delta_t/\omega(t)\to0$ weak-star in $L^{\infty}(1/\omega)$ as $t\to\infty$.
\item
$\text{ran}\,T_{\varphi}\subseteq C_0(1/\omega)$.
\item
$\text{ran}\,T_{\varphi}\subseteq C_0(1/\omega)$ and $T_{\varphi}^*=\overline{D}_{\varphi}$.
\end{enumerate}
\end{proposition}

\noindent{\bf Proof}\quad\\
(a)$\Rightarrow$(b):\quad
This follows because $\delta_t/\omega(t)\to0$ weak-star in $M(\omega)$ as $t\to\infty$.\\
(b)$\Rightarrow$(c):\quad
This follows from Proposition~\ref{pr:tfi}(b) because 
$(T_{\varphi}f)(t)/\omega(t)=\langle f,\overline{D}_{\varphi}\delta_t/\omega(t)\rangle$ 
for $f\in L^1(\omega)$ and $t\in\mathbb R^+$.\\
(c)$\Rightarrow$(d):\quad
For $f,g\in L^1(\omega)$ it follows from (\ref{eq:tfi}) that
$$\langle f,D_{\varphi}g\rangle=\langle T_{\varphi}f,g\rangle$$
since $T_{\varphi}f\in C_0(1/\omega)$ and $g\in L^1(\omega)\le M(\omega)$. Hence $T_{\varphi}^*=D_{\varphi}$ on $L^1(\omega)$.
Let $\mu\in M(\omega)$. By Corollary~\ref{co:bai} we have
$e_n\ast\mu\to\mu$ strongly in $M(\omega)$ as $n\to\infty$, so it follows from Theorem~\ref{th:gr} that $D_{\varphi}(e_n\ast\mu)\to\overline{D}_{\varphi}(\mu)$ weak-star in $L^{\infty}(1/\omega)$ as $n\to\infty$. Moreover, for 
$f\in L^1(\omega)$ we have $e_n\cdot T_{\varphi}f\to T_{\varphi}f$ in $C_0(1/\omega)$ as $n\to\infty$ by Lemma~\ref{le:bai}(a) and thus
$$\langle f,T_{\varphi}^*(e_n\ast\mu)\rangle=\langle T_{\varphi}f,e_n\ast\mu\rangle
=\langle e_n\cdot T_{\varphi}f,\mu\rangle\to\langle T_{\varphi}f,\mu\rangle
=\langle f,T_{\varphi}^*\mu\rangle$$
as $n\to\infty$. Hence 
$D_{\varphi}(e_n\ast\mu)=T_{\varphi}^*(e_n\ast\mu)\to T_{\varphi}^*\mu$ weak-star in $L^{\infty}(1/\omega)$ as $n\to\infty$. It follows that $T_{\varphi}^*\mu=\overline{D}_{\varphi}\mu$ and thus $T_{\varphi}^*=\overline{D}_{\varphi}$ on $M(\omega)$.\\
{[Alternatively, a direct but lengthy calculation shows that $T_{\varphi}^*$ is a derivation. Since $T_{\varphi}^*=D_{\varphi}$ on $L^1(\omega)$, it follows from the uniqueness of the extension from Theorem~\ref{th:gr} that $T_{\varphi}^*=D_{\varphi}$ on $M(\omega)$.]}\\
(d)$\Rightarrow$(a):\quad
Is clear.
{\nopagebreak\hfill\raggedleft$\Box$\bigskip}

We will now show that a certain relatively easily verified condition on $\varphi$ ensures that the equivalent conditions in Proposition~\ref{pr:wkscts} hold. The idea is that if $\varphi/\omega$ is not 
bounded away from zero on large sets, then the definition of $T_{\varphi}f$ can be used to show that $(T_{\varphi}f)(t)/\omega(t)\to0$ as $t\to\infty$ for $f\in L^1(\omega)$.
For $\varphi\in L^{\infty}(1/\omega)$ and $t,\varepsilon\in\mathbb R^+$ we let
$$U_{t,\varepsilon}=\{s\in[t,t+1]\,:\,|\varphi(s)|/\omega(s)\ge\varepsilon\}$$
(defined except for a set of measure zero). 
Also, we denote the Lebesgue measure on $\mathbb R^+$ by $m$. 

\begin{theorem}
\label{th:wkscts}
Let $\varphi\in L^{\infty}(1/\omega)$ and assume that $m(U_{t,\varepsilon})\to0$ as $t\to\infty$ for every $\varepsilon>0$. Then $\overline{D}_{\varphi}$ is weak-star continuous (and consequently the other equivalent conditions in Proposition~\ref{pr:wkscts} also hold).
\end{theorem}

We will need the following lemma.

\begin{lemma}
\label{le:un}
Let $f\in L^1[0,1]$ and let $(U_n)$ be a sequence of measurable sets in $[0,1]$ with $m(U_n)\to0$ as $n\to\infty$. Then $\int_{U_n}f(t)\,dt\to0$ as $n\to\infty$.
\end{lemma}

\noindent{\bf Proof}\quad
It is sufficient to prove that every subsequence $(U_{n_k})$ of $(U_n)$ has a subsequence $(U_{n_{k_j}})$ with 
$\int_{U_{n_{k_j}}}f(t)\,dt\to0$ as $j\to\infty$. We may therefore assume that $\sum_{n=1}^{\infty}m(U_n)<\infty$. Then $m(\cup_{m=n}^{\infty}U_m)\le\sum_{m=n}^{\infty}m(U_m)\to0$ as $n\to\infty$.
Let $f_n=f\cdot1_{U_n}$ for $n\in\mathbb N$. It then follows from 
\cite[Theorem~A, p.\,91]{Ha} that $f_n\to0$ a.e. Consequently 
$\int_{U_n}f(t)\,dt=\int_0^1f_n(t)\,dt\to0$ as $n\to\infty$ by Lebesgue's dominated convergence theorem.
{\nopagebreak\hfill\raggedleft$\Box$\bigskip}

\noindent{\bf Proof of Theorem~\ref{th:wkscts}}\quad
By Proposition~\ref{pr:wkscts} we only need to prove that $\text{ran}\,T_{\varphi}\subseteq C_0(1/\omega)$. (A similar proof can be given to show directly that condition (b) in Proposition~\ref{pr:wkscts} holds.)
We first let $f\in L^1(\omega)$ with $\text{supp}\,f\subseteq[0,1]$. Then
$$|(T_{\varphi}f)(t)|\le\int_0^1|f(s)\varphi(t+s)|\,ds
=\int_t^{t+1}|f(r-t)\varphi(r)|\,dr$$ 
for $t\in\mathbb R^+$. Let $\varepsilon>0$. Then
\begin{align*}
\int_{[t,t+1]\setminus U_{t,\varepsilon}}|f(r-t)\varphi(r)|\,dr
&\le \varepsilon\int_t^{t+1}|f(r-t)|\omega(r)\,dr\\
&\le \varepsilon\int_t^{t+1}|f(r-t)|\omega(r-t)\,dr\,\omega(t)
=\varepsilon\|f\|\omega(t)
\end{align*}
for $t\in\mathbb R^+$. Moreover, 
$$\int_{U_{t,\varepsilon}}|f(r-t)\varphi(r)|\,dr
=\int_{U_{t,\varepsilon}-t}|f(s)\varphi(t+s)|\,ds
\le\|\varphi\|\omega(t)\int_{U_{t,\varepsilon}-t}|f(s)|\omega(s)\,ds$$
for $t\in\mathbb R^+$. It follows from Lemma~\ref{le:un} that there exists $T\in\mathbb R^+$ such that 
$$\int_{U_{t,\varepsilon}-t}|f(s)|\omega(s)\,ds<\varepsilon$$ 
for $t\ge T$. Hence there is a constant $C$ such that 
$|(T_{\varphi}f)(t)|\le C\varepsilon\omega(t)$ for $t\ge T$, so we conclude that $T_{\varphi}f\in C_0(1/\omega)$. 

Next, we let $f\in L^1(\omega)$ with $\text{supp}\,f\subseteq[n,n+1]$ for some $n\in\mathbb N$. Then $f=\delta_n*g$ for some $g\in L^1(\omega)$ with $\text{supp}\,g\subseteq[0,1]$. Also,
$$(T_{\varphi}f)(t)
= \int_n^{n+1}f(s)\,\frac{t}{t+s}\,\varphi(t+s)\,ds
= \int_0^1 g(r)\,\frac{t}{t+r+n}\,\varphi(t+r+n)\,dr,$$
so
$$|(T_{\varphi}f)(t)|
\le\int_0^1|g(r)|\,\frac{t+n}{t+r+n}\,|\varphi(t+r+n)|\,dr
=(T_{|\varphi|}|g|)(t+n)$$
for $t\in\mathbb R^+$. Applying the first part of the proof (the definition of $U_{t,\varepsilon}$ only depends on $|\varphi|$) we get
$$\frac{|(T_{\varphi}f)(t)|}{\omega(t)}
\le\frac{|(T_{|\varphi|}|g|)(t+n)|}{\omega(t+n)}\cdot\omega(n)\to0$$
as $t\to\infty$, so $T_{\varphi}f\in C_0(1/\omega)$. Consequently, $T_{\varphi}f\in C_0(1/\omega)$ for every $f\in L^1(\omega)$ with compact support and thus for all $f\in L^1(\omega)$. 
{\nopagebreak\hfill\raggedleft$\Box$\bigskip}

\begin{corollary}
\label{co:wkscts}
Let $\varphi\in L_0^{\infty}(1/\omega)$. Then $\overline{D}_{\varphi}$ is weak-star continuous and $\overline{D}_{\varphi}=T_{\varphi}^*$.
\end{corollary}

\noindent{\bf Proof}\quad
Let $\varepsilon>0$. There exists $T\in\mathbb R^+$ such that $U_{t,\varepsilon}$ is of measure zero for $t\ge T$. The result thus follows from Theorem~\ref{th:wkscts}.
{\nopagebreak\hfill\raggedleft$\Box$\bigskip}

The following corollary shows a class of functions $\varphi\notin L_0^{\infty}(1/\omega)$ for which $\overline{D}_{\varphi}$ is weak-star continuous.

\begin{corollary}
\label{co:wkscts2}
Let $(\alpha_n)$ be a sequence with $0<\alpha_n<1$ for $n\in\mathbb N$ and $\alpha_n\to0$ for $n\to\infty$. Define $\varphi\in L^{\infty}(1/\omega)$ by the weak-star convergent series $\varphi=\sum_{n=1}^{\infty}1_{[n,n+\alpha_n]}\cdot\omega$. Then $\varphi\notin L_0^{\infty}(1/\omega)$, but $\overline{D}_{\varphi}$ is weak-star continuous.
\end{corollary}

\noindent{\bf Proof}\quad
Let $0<\varepsilon<1$. For $t\in\mathbb R^+$ we let $n=[t]$. Then
$$U_{t,\varepsilon}\subseteq[n,n+\alpha_n]\cup[n+1,n+1+\alpha_{n+1}],$$ 
so $m(U_{t,\varepsilon})\le \alpha_n+\alpha_{n+1}\to0$ as $t\to\infty$. The result thus follows from Theorem~\ref{th:wkscts}.
{\nopagebreak\hfill\raggedleft$\Box$\bigskip}

We do not know whether the condition in Theorem~\ref{th:wkscts} also is necessary for $\overline{D}_{\varphi}$ to be weak-star continuous, but we finish this section by giving a partial result in this direction.

\begin{proposition}
\label{pr:notwkscts}
Suppose that there exists a positive constant $C$ such that $\int_x^{x+1}\omega(y)\,dy\ge C\omega(x)$ for all $x\in\mathbb R^+$. Let $(a_n)$ be a sequence in $\mathbb R^+$ with $a_0\ge1$ and $a_{n+1}\ge a_n+1$ for $n\in\mathbb N$. Define $\varphi\in L^{\infty}(1/\omega)$ by the weak-star convergent series $\varphi=\sum_{n=1}^{\infty}1_{[a_n,a_n+1]}\cdot\omega$. Then $\overline{D}_{\varphi}$ is not weak-star continuous.
\end{proposition}

\noindent{\bf Proof}\quad
For $n\in\mathbb N$ we have 
$$\biggl\langle 1_{[0,1]},\frac{\overline{D}_{\varphi}\delta_{a_n}}{\omega(a_n)}\biggr\rangle
=\int_0^1\frac{a_n}{t+a_n}\cdot\frac{\varphi(t+a_n)}{\omega(a_n)}\,dt
=\int_{a_n}^{a_n+1}\frac{a_n}{s}\cdot\frac{\omega(s)}{\omega(a_n)}\,ds
\ge\frac{C}{2}\,.$$
Hence $\overline{D}_{\varphi}\delta_{a_n}/\omega(a_n)$ does not tend to 0 weak-star in $L^{\infty}(1/\omega)$ as $n\to\infty$. Since $\delta_{a_n}/\omega(a_n)\to0$ weak-star in $M(\omega)$ as $n\to\infty$, this shows that $\overline{D}_{\varphi}$ is not weak-star continuous. 
{\nopagebreak\hfill\raggedleft$\Box$\bigskip}

\begin{corollary}
\label{co:notwkscts}
Suppose that there exists a positive constant $C$ such that $\int_x^{x+1}\omega(y)\,dy\ge C\omega(x)$ for all $x\in\mathbb R^+$. Then $\overline{D}_{\omega}$ is not weak-star continuous.
\end{corollary}

We remark that Corollaries~\ref{co:wkscts}, \ref{co:wkscts2} and Proposition~\ref{pr:notwkscts} can be combined to yield a wider class of functions $\varphi$ for which $\overline{D}_{\varphi}$ is not weak-star continuous. Namely, if 
$\varphi=\varphi_1+\varphi_2$ with $\overline{D}_{\varphi_1}$ weak-star continuous and $\overline{D}_{\varphi_2}$ not weak-star continuous, then $\overline{D}_{\varphi}$ is not weak-star continuous.

We finish the section by showing that $\overline{D}_{\varphi}\mu$ can be represented as a weak-star Bochner integral for 
$\varphi\in L^{\infty}(1/\omega)$ and $\mu\in M(\omega)$. Heuristically we can think of $\overline{D}_{\varphi}\mu$ as
$$(\overline{D}_{\varphi}\mu)(t)
=\int_0^{\infty}\frac{s}{t+s}\,\varphi(t+s)\,d\mu(s)
=\int_0^{\infty}(\overline{D}_{\varphi}\delta_s)(t)\,d\mu(s)
\qquad\text{for $t\in\mathbb R^+$},$$
although the integrals need not be defined. Inspired by this, we will say that
$$\overline{D}_{\varphi}\mu=\int_0^{\infty}\overline{D}_{\varphi}\delta_s\,d\mu(s)$$
as a weak-star Bochner integral in $L^{\infty}(1/\omega)$ if
$$\langle f,\overline{D}_{\varphi}\mu\rangle
=\int_0^{\infty}\langle f,\overline{D}_{\varphi}\delta_s\rangle\,d\mu(s)
\qquad\text{for $f\in L^1(\omega)$}.$$
We remark that the function 
$(T_{\varphi}f)(s)=\langle f,\overline{D}_{\varphi}\delta_s\rangle$ belongs to $C_b(1/\omega)$ by Proposition~\ref{pr:tfi}(b) for $f\in L^1(\omega)$. Hence 
$$\int_0^{\infty}\langle f,\overline{D}_{\varphi}\delta_s\rangle\,d\mu(s)
=\langle T_{\varphi}f,\mu\rangle$$
is well-defined for $f\in L^1(\omega)$ and $\mu\in M(\omega)$. 
Also, it follows from the proof of \cite[Theorem~5.6.24]{Da:Book} that $$D_{\varphi}g=\int_0^{\infty}g(s)\overline{D}_{\varphi}\delta_s\,ds$$ 
as a weak-star integral in $L^{\infty}(1/\omega)$ for every 
$\varphi\in L^{\infty}(1/\omega)$ and $g\in L^1(\omega)$. We will show that this result can be extended to $\overline{D}_{\varphi}\mu$.
(Moreover, in the proof of Proposition~\ref{pr:cpt2} we will see that if $\varphi\in C_b(1/\omega)$ with $\varphi(0)=0$, then 
$D_{\varphi}g=\int_0^{\infty}g(s)\overline{D}_{\varphi}\delta_s\,ds$ actually exists as a ``proper'' Bochner integral for $g\in L^1(\omega)$.)

\begin{proposition}
\label{pr:Bochner}
Let $\varphi\in L^{\infty}(1/\omega)$ and $\mu\in M(\omega)$. Then $$\overline{D}_{\varphi}\mu=\int_0^{\infty}\overline{D}_{\varphi}\delta_s\,d\mu(s)$$
as a weak-star Bochner integral in $L^{\infty}(1/\omega)$.
\end{proposition}

\noindent{\bf Proof}\quad
Let $(e_n)$ be the bounded approximate identity from Lemma~\ref{le:bai}. By Corollary~\ref{co:bai} we have $e_n\ast\mu\to\mu$ strongly in $M(\omega)$ as $n\to\infty$, so it follows from Theorem~\ref{th:gr} that  $D_{\varphi}(e_n\ast\mu)\to\overline{D}_{\varphi}(\mu)$ weak-star in $L^{\infty}(1/\omega)$ as $n\to\infty$. Since $e_n\ast\mu\in L^1(\omega)$ we have
$$D_{\varphi}(e_n\ast\mu)
=\int_0^{\infty}(e_n\ast\mu)(s)\overline{D}_{\varphi}\delta_s\,ds$$ 
as a weak-star integral in $L^{\infty}(1/\omega)$ for $n\in\mathbb N$.
By using Corollary~\ref{co:bai} we thus have
\begin{align*}
\langle f,\overline{D}_{\varphi}\mu\rangle
&=\lim_{n\to\infty}\langle f,D_{\varphi}(e_n\ast\mu)\rangle\\
&=\lim_{n\to\infty}\int_0^{\infty}
\langle f,\overline{D}_{\varphi}\delta_s\rangle(e_n\ast\mu)(s)\,ds\\
&=\lim_{n\to\infty}\langle T_{\varphi}f,e_n\ast\mu\rangle
=\langle T_{\varphi}f,\mu\rangle
=\int_0^{\infty}\langle f,\overline{D}_{\varphi}\delta_s\rangle\,d\mu(s)
\end{align*}
as required.
{\nopagebreak\hfill\raggedleft$\Box$\bigskip}

\section{Range and continuity properties}
\label{sec:various}

For $\varphi\in L^{\infty}(1/\omega)$ we saw in (the proof of) Proposition~\ref{pr:tfi}(a) that $\overline{D}_{\varphi}\delta_s\to0$ weak-star in $L^{\infty}(1/\omega)$ as $s\to0$. Similarly, one of the equivalent conditions in Proposition~\ref{pr:wkscts} is that $\overline{D}_{\varphi}\delta_s/\omega(s)\to0$ weak-star in $L^{\infty}(1/\omega)$ as $s\to\infty$. We will now see that by strengthening the conditions on $\varphi$, we can obtain norm convergence in both cases.

\begin{proposition}
\label{pr:lioom}
Let $\varphi\in L^{\infty}(1/\omega)$. 
\begin{enumerate}[(a)]
\item 
$\varphi(s)\to0$ as $s\to0$ if and only if $\overline{D}_{\varphi}\delta_s\to0$ in $L^{\infty}(1/\omega)$ as $s\to0$. 
\item
$\varphi\in L_0^{\infty}(1/\omega)$ if and only if $\overline{D}_{\varphi}\delta_s/\omega(s)\to0$ in $L^{\infty}(1/\omega)$ as $s\to\infty$. 
\end{enumerate}
\end{proposition}

\noindent{\bf Proof}\quad
(a):\quad
Assume that $\varphi(s)\to0$ as $s\to0$. We have
$$\|\overline{D}_{\varphi}\delta_s\|=\text{ess\,sup}_{t\in\mathbb R^+}\,\frac{s}{t+s}\cdot
\frac{|\varphi(t+s)|}{\omega(t)}\,.$$
Given $\varepsilon>0$ we choose $0<S<1$ such that 
$\text{ess\,sup}_{s\le S}\,|\varphi(s)|<\varepsilon$. For $s\le S/2$ we then have
$$\text{ess\,sup}_{t\le S/2}\,\frac{s}{t+s}\cdot
\frac{|\varphi(t+s)|}{\omega(t)}
\le\frac{\varepsilon}{\inf_{t\in[0,1]}\omega(t)}\,.$$
Also, since $|\varphi(t+s)|\le\|\varphi\|\omega(t)\omega(s)$ for $t,s\in\mathbb R^+$, we have
$$\text{ess\,sup}_{t\ge S/2}\,\frac{s}{t+s}\cdot
\frac{|\varphi(t+s)|}{\omega(t)}
\le\frac{2s}{S}\|\varphi\|\omega(s),$$
and it follows that $\overline{D}_{\varphi}\delta_s\to0$ in $L^{\infty}(1/\omega)$ as $s\to0$. 

Conversely, assume that $\overline{D}_{\varphi}\delta_s\to0$ in $L^{\infty}(1/\omega)$ as $s\to0$. Given $\varepsilon>0$ we choose $0<S<1$ such that $\|\overline{D}_{\varphi}\delta_s\|<\varepsilon$ for $s\le S$. For $s\le S$ we then have
$$\varepsilon>\|\overline{D}_{\varphi}\delta_s\|
\ge\frac{\text{ess\,sup}_{t\le s}\,|\varphi(t+s)|}{2\sup_{t\in[0,1]}\omega(t)}
=C\text{ess\,sup}_{s\le t\le 2s}\,|\varphi(t)|.$$
Hence $\text{ess\,sup}_{0\le t\le 2S}\,|\varphi(t)|\le\varepsilon/C$, so 
$\varphi(s)\to0$ as $s\to0$.

(b):\quad
Assume that $\varphi\in L_0^{\infty}(1/\omega)$. Then 
\begin{align*}
\frac{\|\overline{D}_{\varphi}\delta_s\|}{\omega(s)}
&=\text{ess\,sup}_{t\in\mathbb R^+}\,\frac{s}{t+s}\cdot
\frac{|\varphi(t+s)|}{\omega(t)\omega(s)}\\
&\le\text{ess\,sup}_{t\in\mathbb R^+}\,\frac{|\varphi(t+s)|}{\omega(t+s)}
=\text{ess\,sup}_{t\ge s}\,\frac{|\varphi(t)|}{\omega(t)}\to0
\end{align*}
as $s\to\infty$. 

Conversely, assume that $\varphi\notin L_0^{\infty}(1/\omega)$. Then there exists $\varepsilon>0$ such that 
$$\text{ess\,sup}_{t\ge T}\,|\varphi(t)|/\omega(t)\ge\varepsilon$$ 
for every $T\in\mathbb R^+$. Let $k\in\mathbb N$.  There exists a measurable set $U_k\subseteq[k,\infty)$ with $m(U_k)>0$ such that $|\varphi(t)|/\omega(t)\ge\varepsilon$ a.e.\ on $U_k$. The metric density of $U_k$ at a point $s\in U_k$, that is $\lim_{r\to0}m(U_k\cap(s-r,s+r))/(2r)$, exists and equals 1 for almost every $s\in U_k$ (\cite[7.12]{Ru}). Let 
$s_k\in U_k$ be such a point and let $V_{kr}=U_k\cap[s_k,s_k+r)$ for $r>0$. Then $m(V_{kr})>0$ for every $r>0$. Hence
\begin{align*}
\|\overline{D}_{\varphi}\delta_{s_k}\|
&\ge\text{ess\,sup}_{t\in[0,1]}\,\frac{s_k}{t+s_k}\cdot
\frac{|\varphi(t+s_k)|}{\omega(t)}\\
&\ge\frac{1}{2\inf_{t\in[0,1]}\omega(t)}\,
\text{ess\,sup}_{t\in[0,1]}\,|\varphi(t+s_k)|
\ge C\text{ess\,sup}_{t\in V_{kr}}\,|\varphi(t)|
\end{align*}
for some constant $C>0$ and every $0<r<1$. Since $V_{kr}\subseteq U_k$ we thus have
$$\|\overline{D}_{\varphi}\delta_{s_k}\|\ge C\varepsilon\sup_{t\in V_{kr}}\omega(t)$$
for every $0<r<1$. Letting $r\to0$ we thus obtain 
$$\|\overline{D}_{\varphi}\delta_{s_k}\|\ge C\varepsilon\omega(s_k).$$
Since $s_k\ge k$ this shows that we do not have $\overline{D}_{\varphi}\delta_s/\omega(s)\to0$ in $L^{\infty}(1/\omega)$ as $s\to\infty$. 
{\nopagebreak\hfill\raggedleft$\Box$\bigskip}

We now aim to prove that $\text{ran}\,D_{\varphi}\subseteq C_0(1/\omega)$ for any $\varphi\in L^{\infty}(1/\omega)$, and that $\text{ran}\,\overline{D}_{\varphi}\subseteq C_0(1/\omega)$ if 
$\varphi\in C_0(1/\omega)$ with $\varphi(0)=0$. Let 
$\varphi\in L^{\infty}(1/\omega)$. Since 
$(D_{\varphi}f)(t)=\int_0^{\infty}\frac{s}{t+s}\,f(s)\varphi(t+s)\,ds$
for $f\in L^1(\omega)$ and $t\in\mathbb R^+$, we choose to define
$$\psi_t(s)=(\overline{D}_{\varphi}\delta_s)(t)=\frac{s}{t+s}\,\varphi(t+s)\qquad
\text{for } t,s\in\mathbb R^+.$$
Then we can express $D_{\varphi}f$ by
$$(D_{\varphi}f)(t)=\langle f,\psi_t\rangle\qquad
\text{for $f\in L^1(\omega)$ and $t\in\mathbb R^+$}$$
(once we have verified that $\psi_t\in L^{\infty}(1/\omega)$ for $t\in\mathbb R^+$).
We begin by establishing som properties of $\psi_t$.

\begin{lemma}
\label{le:psi}
\ 
\begin{enumerate}[(a)]
\item 
Let $\varphi\in L^{\infty}(1/\omega)$. For $t\in\mathbb R^+$ we have $\psi_t\in L^{\infty}(1/\omega)$ with $\|\psi_t\|\le\|\varphi\|\omega(t)$. Moreover, $(\psi_t)$ is weak-star continuous in $L^{\infty}(1/\omega)$ for $t\in\mathbb R^+$.
\item
Let $\varphi\in C_0(1/\omega)$ with $\varphi(0)=0$. For $t\in\mathbb R^+$ we have $\psi_t\in C_0(1/\omega)$. Moreover, $(\psi_t)$ is continuous in $C_0(1/\omega)$ for $t\in\mathbb R^+$ and $\psi_t/\omega(t)\to0$ in $C_0(1/\omega)$ as $t\to\infty$. 
\end{enumerate}
\end{lemma}

\noindent{\bf Proof}\quad
(a):\quad
Let $t\in\mathbb R^+$. We have $$\frac{|\psi_t(s)|}{\omega(s)}\le\frac{|\varphi(t+s)|}{\omega(t+s)}\,\omega(t)
\le\|\varphi\|\omega(t)\qquad\text{for all }s\in\mathbb R^+,$$ 
so $\psi_t\in L^{\infty}(1/\omega)$ with $\|\psi_t\|\le\|\varphi\|\omega(t)$. Also, translation is weak-star continuous in $L^{\infty}(1/\omega)$, so $(\psi_t)$ is weak-star continuous in $L^{\infty}(1/\omega)$ for $t>0$. Also, for $f\in L^1(\omega)$ we have
\begin{align*}
\langle f,\psi_t-\psi_0\rangle
&=\int_0^{\infty}f(s)\left(\frac{s}{t+s}\,\varphi(t+s)-\varphi(s)\right)\,ds\\
&=\int_0^{\infty}f(s)(\varphi(t+s)-\varphi(s))\,ds
-\int_0^{\infty}f(s)\,\frac{t}{t+s}\,\varphi(t+s)\,ds.
\end{align*}
As $t\to0$ the first term tends to 0 since translation is weak-star continuous in $L^{\infty}(1/\omega)$, whereas the second tends to 0 by Lebesgue's dominated convergence theorem since $|\varphi(t+s)|\le\|\varphi\|\omega(t)\omega(s)$ for $t,s\in\mathbb R^+$. Hence $(\psi_t)$ is also weak-star continuous in $L^{\infty}(1/\omega)$ at $t=0$.

(b):\quad
Clearly $\psi_t\in C_0(1/\omega)$ for $t\in\mathbb R^+$, and since translation is continuous in $C_0(1/\omega)$, it follows that $(\psi_t)$ is continuous in $C_0(1/\omega)$ for $t>0$. We will now prove that $(\psi_t)$ is also continuous in $C_0(1/\omega)$ at $t=0$. For $t\in\mathbb R^+$ we have
$$\|\psi_t-\psi_0\|=\text{ess\,sup}_{s\in\mathbb R^+}\,
\frac{\left|\frac{s}{t+s}\,\varphi(t+s)-\varphi(s)\right|}{\omega(s)}\,.$$
Given $\varepsilon>0$, we choose $S_1,S_2>0$ such that $|\varphi(s)|/\omega(s)\le\varepsilon$ if $s\le S_1$ or $s\ge S_2$. For 
$t\in\mathbb R^+$ we then have
$$\text{ess\,sup}_{s\ge S_2}\, \frac{\left|\frac{s}{t+s}\,\varphi(t+s)-\varphi(s)\right|}{\omega(s)}
\le\text{ess\,sup}_{s\ge S_2} \left(\frac{\omega(t+s)}{\omega(s)}+1\right)\varepsilon
\le(\omega(t)+1)\varepsilon.$$
Also, for $t\in\mathbb R^+$ we  have
\begin{gather*}
\text{ess\,sup}_{s\le S_2}\,
\frac{\left|\frac{s}{t+s}\,\varphi(t+s)-\varphi(s)\right|}{\omega(s)}
\le\text{ess\,sup}_{s\le S_2}\, \frac{s}{t+s}\cdot\frac{|\varphi(t+s)-\varphi(s)|}{\omega(s)}\\
+\text{ess\,sup}_{s\le S_1}\,\frac{t}{t+s}\cdot\frac{|\varphi(s)|}{\omega(s)}
+\text{ess\,sup}_{S_1\le s\le S_2}\, \frac{t}{t+s}\cdot\frac{|\varphi(s)|}{\omega(s)}\,.
\end{gather*}
The first term tends to 0 as $t\to0$ since translation is continuous in $C_0(1/\omega)$; whereas the second is bounded by $\epsilon$ and the third by $Ct$ for some constant $C$. Together these estimates show that $\psi_t\to\psi_0$ in $C_0(1/\omega)$ as $t\to0$ as required.
Moreover, 
$$\frac{\|\psi_t\|}{\omega(t)}
\le\text{ess\,sup}_{s\in\mathbb R^+}\,\frac{|\varphi(t+s)|}{\omega(t)\omega(s)}
\le\text{ess\,sup}_{s\in\mathbb R^+}\,\frac{|\varphi(t+s)|}{\omega(t+s)}
=\text{ess\,sup}_{s\ge t}\,\frac{|\varphi(s)|}{\omega(s)}\to0$$
as $t\to\infty$. 
{\nopagebreak\hfill\raggedleft$\Box$\bigskip}

\begin{corollary}
\label{co:range}
\ 
\begin{enumerate}[(a)]
\item 
Let $\varphi\in L^{\infty}(1/\omega)$. Then 
$(D_{\varphi}f)(t)=\langle f,\psi_t\rangle$ for $f\in L^1(\omega)$ and $t\in\mathbb R^+$. Moreover, $\text{ran}\,D_{\varphi}\subseteq C_0(1/\omega)$.
\item
Let $\varphi\in C_0(1/\omega)$ with $\varphi(0)=0$. Then 
$(\overline{D}_{\varphi}\mu)(t)=\langle \psi_t,\mu\rangle$ for $\mu\in M(\omega)$ and $t\in\mathbb R^+$. Moreover, $\text{ran}\,\overline{D}_{\varphi}\subseteq C_0(1/\omega)$. 
\end{enumerate}
\end{corollary}

\noindent{\bf Proof}\quad
\quad
(a):\quad 
Let $f\in L^1(\omega)$. It follows from Lemma~\ref{le:psi}(a) that 
$(D_{\varphi}f)(t)=\langle f,\psi_t\rangle$ for $t\in\mathbb R^+$. From this it follows  that $D_{\varphi}f$ is continuous on $\mathbb R^+$ and that 
$|(D_{\varphi}f)(t)|\le\|f\|\cdot\|\psi_t\|\le\|f\|\cdot\|\varphi\|\omega(t)$. Hence $D_{\varphi}f\in C_b(1/\omega)$. Let $\varepsilon>0$ and choose $S>0$ such that 
$\int_S^{\infty}|f(s)|\omega(s)\,ds<\varepsilon$. For $t\ge S/\varepsilon$ we then have
\begin{align*}
\frac{|(D_{\varphi}f)(t)|}{\omega(t)}
&\le\|\varphi\|
\left(\int_0^S|f(s)|\,\frac{s}{t+s}\,\frac{\omega(t+s)}{\omega(t)}\,ds
+\int_S^{\infty}|f(s)|\,\frac{\omega(t+s)}{\omega(t)}\,ds\right)\\
&\le\|\varphi\|
\left(\int_0^S\varepsilon|f(s)|\omega(s)\,ds
+\int_S^{\infty}|f(s)|\omega(s)\,ds\right)\le C\varepsilon
\end{align*}
for some constant $C$, so we conclude that $D_{\varphi}f\in C_0(1/\omega)$.

(b):\quad
Let $\mu\in M(\omega)$ and define
$$k(t)=\int_0^{\infty}\frac{s}{t+s}\,\varphi(t+s)\,d\mu(s)
=\langle\psi_t,\mu\rangle\qquad\text{for }t\in\mathbb R^+.$$
Then $k\in C_0(1/\omega)$ by Lemma~\ref{le:psi}(b). By Corollary~\ref{co:wkscts} we have $\overline{D}_{\varphi}=T_{\varphi}^*$, so for $f\in L^1(\omega)$ we have 
\begin{align*}
\langle f,\overline{D}_{\varphi}\mu\rangle
&= \langle T_{\varphi}f,\mu\rangle
=\int_0^{\infty}(T_{\varphi}f)(s)\,d\mu(s)
=\int_0^{\infty}\int_0^{\infty}f(t)\,\frac{s}{t+s}\,\varphi(t+s)\,dt\,d\mu(s)\\
&= \int_0^{\infty}\int_0^{\infty}\frac{s}{t+s}\,\varphi(t+s)\,d\mu(s)\,f(t)\,dt
=\int_0^{\infty}k(t)f(t)\,dt=\langle f,k\rangle.
\end{align*}
Hence $\overline{D}_{\varphi}\mu=k\in C_0(1/\omega)$ and the conclusions follow.
{\nopagebreak\hfill\raggedleft$\Box$\bigskip}

We finish the section with the following result, which will be used in the next section.

\begin{proposition}
\label{pr:cbom}
Let $\varphi\in C_b(1/\omega)$ with $\varphi(0)=0$. Then $(\overline{D}_{\varphi}\delta_s)$ is continuous in $C_0(1/\omega)$ for $s\in\mathbb R^+$.
\end{proposition}

\noindent{\bf Proof}\quad
Clearly $\overline{D}_{\varphi}\delta_s\in C_0(1/\omega)$ for $s\in\mathbb R^+$. Let $s_0>0$. For $x>-s_0$ and $t\in\mathbb R^+$ we have
\begin{align*}
(\overline{D}_{\varphi}\delta_{s_0+x})(t)
&=\frac{s_0+x}{t+s_0+x}\,\varphi(t+s_0+x)\\
&=\frac{s_0+x}{s_0}\,(\overline{D}_{\varphi}\delta_{s_0})(t+x)
=\frac{s_0+x}{s_0}\,(\delta_{-x}*\overline{D}_{\varphi}\delta_{s_0})(t).
\end{align*}
Since translation is continuous in $C_0(1/\omega)$ we have $\delta_{-x}*\overline{D}_{\varphi}\delta_{s_0}\to\overline{D}_{\varphi}\delta_{s_0}$ and thus $\overline{D}_{\varphi}\delta_{s_0+x}\to\overline{D}_{\varphi}\delta_{s_0}$ in $C_0(1/\omega)$ as $x\to0$. Hence $(\overline{D}_{\varphi}\delta_s)$ is continuous in $C_0(1/\omega)$ at $s_0$. 
Finally, by Proposition~\ref{pr:lioom}(a) we have $\overline{D}_{\varphi}\delta_s\to\overline{D}_{\varphi}\delta_0=0$ in $C_0(1/\omega)$ as $s\to0$, so $(\overline{D}_{\varphi}\delta_s)$ is also continuous in $C_0(1/\omega)$ at $s=0$. 
{\nopagebreak\hfill\raggedleft$\Box$\bigskip}

\section{Compactness}
\label{sec:compact}

In this section we study compactness of the operators $D_{\varphi}$ and  $\overline{D}_{\varphi}$. The main result of the section is Theorem~\ref{th:cpt}, which states that for $\varphi\in C_0(1/\omega)$ the operator $\overline{D}_{\varphi}$ is compact if and only if $\varphi(0)=0$.
We start with some results which shows why the condition $\varphi(0)=0$ as well as the continuity of $\varphi$ seem to be close to necessary.

\begin{proposition}
\label{pr:cpt}
Let $\varphi\in L^{\infty}(1/\omega)$ and assume that $\overline{D}_{\varphi}$ is compact. Then $\varphi(s)\to0$ as $s\to0$. 
\end{proposition}

\noindent{\bf Proof}\quad
It follows from Proposition~\ref{pr:tfi}(a) that $\overline{D}_{\varphi}\delta_s\to0$ weak-star in $L^{\infty}(1/\omega)$ as $s\to0$. Since $\overline{D}_{\varphi}$ is compact, we then also have $\overline{D}_{\varphi}\delta_s\to0$ in norm in $L^{\infty}(1/\omega)$ as $s\to0$. Hence  
$\varphi(s)\to0$ as $s\to0$ by Proposition~\ref{pr:lioom}(a).
{\nopagebreak\hfill\raggedleft$\Box$\bigskip}

\begin{theorem}
\label{th:notcpt2}
Let $\varphi\in L^{\infty}(1/\omega)$ be real-valued and assume that there exist $t_0,\delta>0$ such that
$$\text{ess\,inf}_{t\in(t_0-\delta,t_0)}\,\varphi(t)
>\text{ess\,sup}_{t\in(t_0,t_0+\delta)}\,\varphi(t).$$
Then $D_{\varphi}$ and $\overline{D}_{\varphi}$ are not compact.
\end{theorem}

\noindent{\bf Proof}\quad
Choose $\alpha\in\mathbb R$ and $\varepsilon>0$ such that
$$\varphi(t)\ge\alpha+\varepsilon\quad\text{a.e.\ on }(t_0-\delta,t_0)
\qquad\text{and}\qquad 
\varphi(t)\le\alpha-\varepsilon\quad\text{a.e.\ on }(t_0,t_0+\delta).$$
Assume that $D_{\varphi}$ is compact and let $f_n=n\cdot1_{[t_0-1/n,t_0]}$ for $n\in\mathbb N$ with $1/n\le t_0$. Then $(D_{\varphi}f_n)$ has a norm-convergent subsequence $(D_{\varphi}f_{n_k})$ with limit $h\in L^{\infty}(1/\omega)$. For $n\in\mathbb N$ with 
$1/n\le t_0$ and $t\in\mathbb R^+$  we have
$$(D_{\varphi}f_n)(t)=n\int_{t_0-1/n}^{t_0}\frac{s}{t+s}\,\varphi(t+s)\,ds,$$
so $(D_{\varphi}f_n)(t)\le\alpha-\varepsilon$ for $1/n\le t\le\delta$. Hence $h(t)\le\alpha-\varepsilon$ a.e.\ on $[0,\delta]$. 
For $n\in\mathbb N$ with $1/n\le t_0$ and $t\le 1/n$ we have 
$$(D_{\varphi}f_n)(t)
=n\int_{t_0-1/n}^{t_0-t}\frac{s}{t+s}\,\varphi(t+s)\,ds
+n\int_{t_0-t}^{t_0}\frac{s}{t+s}\,\varphi(t+s)\,ds
=A_n(t)+B_n(t)$$
with obvious notation. We have
$$A_n(t)\ge n(\tfrac{1}{n}-t)\,\frac{t_0-1/n}{t+t_0-1/n}(\alpha+\varepsilon)
\ge (1-nt)\,\frac{t_0-1/n}{t_0}(\alpha+\varepsilon),$$
so there exists $N_1\in\mathbb N$ such that 
$A_n(t)\ge\alpha+\varepsilon/2$ for $n\ge N_1$ and $t\le 1/n^2$. Also,
$$|B_n(t)|\le nt\|\varphi\|\sup_{s\le t_0+\delta}\omega(s),$$
so there exists $N_2\in\mathbb N$ such that 
$|B_n(t)|\le\varepsilon/2$ for $n\ge N_2$ and $t\le 1/n^2$. Hence 
$(D_{\varphi}f_n)(t)\ge\alpha$ for $n\ge\max\{N_1,N_2\}$ and $t\le 1/n^2$.
Hence $\|D_{\varphi}f_n-h\|\ge\varepsilon$ for $n\ge\max\{N_1,N_2\}$ which contradicts $D_{\varphi}f_{n_k}\to h$ in $L^{\infty}(1/\omega)$ a $k\to\infty$. Hence $D_{\varphi}$ is not compact and as a consequence $\overline{D}_{\varphi}$ is not compact either. 
{\nopagebreak\hfill\raggedleft$\Box$\bigskip}

The following corollary implies that simple functions like $\varphi=1_{[0,1]}$ do not generate compact derivations.

\begin{corollary}
\label{co:notcpt2}
Let $\varphi\in L^{\infty}(1/\omega)$ be real-valued and assume that there exists $t_0>0$ such that the limits
$$\lim_{t\to(t_0)_-}\varphi(t)\qquad\text{and}\qquad\lim_{t\to(t_0)_+}\varphi(t)$$ exist and are different. Then $D_{\varphi}$ and $\overline{D}_{\varphi}$ are not compact.
\end{corollary}

\noindent{\bf Proof}\quad
The result follows directly from Theorem~\ref{th:notcpt2} if 
$\lim_{t\to(t_0)_-}\varphi(t)>\lim_{t\to(t_0)_+}\varphi(t)$. If the opposite inequality holds, then the result follows by considering $-\varphi$.
{\nopagebreak\hfill\raggedleft$\Box$\bigskip}

Because of the results above, we will focus on $\varphi\in C_b(1/\omega)$ with $\varphi(0)=0$ in the rest of the paper.
For $\varphi\in C_0(1/\omega)$ we have the following characterisation of compact $\overline{D}_{\varphi}$.

\begin{theorem}
\label{th:cpt}
Let $\varphi\in C_0(1/\omega)$. Then $\overline{D}_{\varphi}$ is compact if and only if $\varphi(0)=0$.
\end{theorem}

\noindent{\bf Proof}\quad
If $\overline{D}_{\varphi}$ is compact, then $\varphi(0)=0$ by Proposition~\ref{pr:cpt} .

Conversely, assume that $\varphi(0)=0$ and let $(\mu_n)$ be a bounded sequence in $M(\omega)$. By passing to subsequences we may assume that there exist $\mu\in M(\omega)$ and $h\in L^{\infty}(1/\omega)$ such that
$$\mu_n\to\mu\quad\text{weak-star in $M(\omega)$\quad and\quad
$\overline{D}_{\varphi}\mu_n\to h$ weak-star in $L^{\infty}(1/\omega)$}$$
as $n\to\infty$. By Corollary~\ref{co:wkscts} we have $\overline{D}_{\varphi}=T_{\varphi}^*$. For $f\in L^1(\omega)$ we thus have 
$$\langle f,h\rangle
=\lim_{n\to\infty}\langle f,\overline{D}_{\varphi}\mu_n\rangle
=\lim_{n\to\infty}\langle T_{\varphi}f,\mu_n\rangle=\langle T_{\varphi}f,\mu\rangle
=\langle f,\overline{D}_{\varphi}\mu\rangle,$$
so we deduce that $h=\overline{D}_{\varphi}\mu$. By Corollary~\ref{co:range}(b) we have
$$(\overline{D}_{\varphi}\mu_n-h)(t)=(\overline{D}_{\varphi}(\mu_n-\mu))(t)=\langle \psi_t,\mu_n-\mu\rangle.$$
Also, $(\psi_t)$ is continuous in $C_0(1/\omega)$ and $\psi_t/\omega(t)\to0$ in $C_0(1/\omega)$ as $t\to\infty$ by Lemma~\ref{le:psi}(b), so $\{\psi_t/\omega(t):t\in\mathbb R^+\}$ is totally bounded in $C_0(1/\omega)$.
Let $\varepsilon>0$. There exist $t_1,\ldots,t_M\in\mathbb R^+$ such that for every $t\in\mathbb R^+$ there exists $m\in\{1,\ldots,M\}$ with 
$\|\psi_t/\omega(t)-\psi_{t_m}/\omega(t_m)\|<\varepsilon$. Choose $N\in\mathbb N$ such that $|\langle\psi_{t_m}/\omega(t_m),\mu_n-\mu\rangle|<\varepsilon$ for $m=1,\ldots,M$ and $n\ge N$. For $t\in\mathbb R^+$ and $n\ge N$ we thus have
\begin{align*}
\frac{|(\overline{D}_{\varphi}\mu_n-h)(t)|}{\omega(t)}
&= \left|\biggl\langle\frac{\psi_t}{\omega(t)},\mu_n-\mu\biggr\rangle\right|\\
&\le\left|\biggl\langle\frac{\psi_{t_m}}{\omega(t_m)},
\mu_n-\mu\biggr\rangle\right|
+\left|\biggl\langle\frac{\psi_t}{\omega(t)}-\frac{\psi_{t_m}}{\omega(t_m)},
\mu_n-\mu\biggr\rangle\right|\\
&< (1+\sup_{n\in\mathbb N}\|\mu_n-\mu\|)\varepsilon.
\end{align*}
Hence $\overline{D}_{\varphi}\mu_n\to h$ in $C_0(1/\omega)$ as $n\to\infty$, and we conclude that $\overline{D}_{\varphi}$ is compact.
{\nopagebreak\hfill\raggedleft$\Box$\bigskip}

The next few results will show the existence of $\varphi\notin C_0(1/\omega)$ for which $D_{\varphi}$ is compact. We do not know whether the approach can be extented to show that $\overline{D}_{\varphi}$ is compact, and more generally we do not know whether $\overline{D}_{\varphi}$ is necessarily compact if $D_{\varphi}$ is compact. The idea in the following result is to represent the derivation $D_{\varphi}$ by Bochner integrals and then use a pre-compactness argument.

\begin{proposition}
\label{pr:cpt2}
Let $\varphi\in C_b(1/\omega)$. Assume that $\varphi(0)=0$ and that $\overline{D}_{\varphi}\delta_s/\omega(s)$ has a limit in $L^{\infty}(1/\omega)$ as $s\to\infty$. Then $D_{\varphi}$ is compact.
\end{proposition}

\noindent{\bf Proof}\quad
Let $f\in L^1(\omega)$. We observe that 
$$(D_{\varphi}f)(t)=\int_0^{\infty}f(s)(\overline{D}_{\varphi}\delta_s)(t)\,ds
\qquad\text{for }t\in\mathbb R^+.$$
Moreover, $(\overline{D}_{\varphi}\delta_s)$ is continuous in $C_0(1/\omega)$ for $s\in\mathbb R^+$ by Proposition~\ref{pr:cbom}, so 
$$D_{\varphi}f=\int_0^{\infty}f(s)\overline{D}_{\varphi}\delta_s\,ds
=\int_0^{\infty}f(s)\omega(s)\,\frac{\overline{D}_{\varphi}\delta_s}{\omega(s)}\,ds$$
exists as a Bochner integral (see \cite[Theorem~3.7.4]{Hi-Ph}). Also, $\overline{D}_{\varphi}\delta_s\to0$ in $L^{\infty}(1/\omega)$ as $s\to0$ by Proposition~\ref{pr:lioom}(a). Hence the map $s\mapsto\overline{D}_{\varphi}\delta_s/\omega(s)$ extends to a continuous map from the one point compactification $[0,\infty]$ of $\mathbb R^+$ into $C_0(1/\omega)$, so we deduce that $\{\overline{D}_{\varphi}\delta_s/\omega(s):s\in\mathbb R^+\}\cup
\{\lim_{s\to\infty}\overline{D}_{\varphi}\delta_s/\omega(s)\}$ is compact in $C_0(1/\omega)$. It thus follows from \cite[VI.8.11]{Du-Sc:I} 
(see also \cite[Proof of Theorem~2.2]{Ba-Da:Norms}) that $D_{\varphi}$ is compact.
{\nopagebreak\hfill\raggedleft$\Box$\bigskip}

For $\varphi\in C_0(1/\omega)$ with $\varphi(0)=0$ we have $\overline{D}_{\varphi}\delta_s/\omega(s)\to0$ in $C_0(1/\omega)$ as $s\to\infty$ by Proposition~\ref{pr:lioom}(b), so $D_{\varphi}$ is compact by Proposition~\ref{pr:cpt2}. We can therefore recapture the conclusion about $D_{\varphi}$ from Theorem~\ref{th:cpt}.

We will now use Proposition~\ref{pr:cpt2} to obtain concrete examples of functions $\varphi\notin C_0(1/\omega)$ which generate compact derivations $D_{\varphi}$.

\begin{proposition}
\label{pr:cpt3}
Assume that $\omega(s)\ge1$ for every $s\in\mathbb R^+$, $\omega(s)\to\infty$ as $s\to\infty$ and
$$\sup_{t\in\mathbb R^+}\frac{|\omega(t+s)-\omega(s)|}{\omega(t)\omega(s)}\to0
\qquad\text{as }s\to\infty.$$ 
Let $\varphi=\omega-1$. Then $D_{\varphi}$ is compact, whereas $\varphi\notin C_0(1/\omega)$.
\end{proposition}

\noindent
{\bf Remark}\quad 
If we let $\rho_s(t)=\omega(t+s)-\omega(s)$ for $t,s\in\mathbb R^+$, then the last assumption can be restated as $\rho_s/\omega(s)\to0$ in $L^{\infty}(1/\omega)$ as $s\to\infty$.  
\bigskip

\noindent{\bf Proof}\quad
By Proposition~\ref{pr:cpt2} the result will follow if we can prove that $\overline{D}_{\varphi}\delta_s/\omega(s)\to1$ in $L^{\infty}(1/\omega)$ as $s\to\infty$. We have
$$\frac{\overline{D}_{\varphi}\delta_s}{\omega(s)}-1
=\frac{\frac{s}{t+s}(\omega(t+s)-1)}{\omega(s)}-1
=\frac{s\omega(t+s)-s-(t+s)\omega(s)}{(t+s)\omega(s)}\,,$$
so
\begin{align}
\label{eq:Dfi}
\left\|\frac{\overline{D}_{\varphi}\delta_s}{\omega(s)}-1\right\|
&=\sup_{t\in\mathbb R^+}\frac{|s(\omega(t+s)-\omega(s))-s-t\omega(s)|}
{(t+s)\omega(t)\omega(s)}
\nonumber\\
&\le\sup_{t\in\mathbb R^+}\frac{|\omega(t+s)-\omega(s)|}{\omega(t)\omega(s)}
+\frac{1}{\omega(s)}+\sup_{t\in\mathbb R^+}\frac{t}{(t+s)\omega(t)}\,.
\end{align}
Given $\varepsilon>0$ we choose $T>0$ such that $\omega(t)\ge1/\varepsilon$ for $t\ge T$. Then 
$$\sup_{t\ge T}\frac{t}{(t+s)\omega(t)}\le\varepsilon$$ 
for every $s\in\mathbb R^+$. Also, for $s\ge T/\varepsilon$ we have 
$$\sup_{t\le T}\frac{t}{(t+s)\omega(t)}
\le\frac{T}{T+T/\varepsilon}<\varepsilon.$$ 
Hence the third term in (\ref{eq:Dfi}) tends to 0 as $s\to\infty$. 
Moreover, the first and second term tend to 0 as $s\to\infty$ by assumption, so we conclude that $\overline{D}_{\varphi}\delta_s/\omega(s)\to1$ in $L^{\infty}(1/\omega)$ as $s\to\infty$. 
{\nopagebreak\hfill\raggedleft$\Box$\bigskip}

\begin{corollary}
\label{co:cpt2}
Let $\alpha>0,\ \omega(t)=(1+t)^{\alpha}\ (t\in\mathbb R^+)$ and let $\varphi=\omega-1$. Then $D_{\varphi}$ is compact, whereas $\varphi\notin C_0(1/\omega)$.
\end{corollary}

\noindent{\bf Proof}\quad
We clearly have $\varphi\notin C_0(1/\omega)$.
First, assume that $\alpha\ge1$. Then
\begin{align*}
0\le\omega(t+s)-\omega(s)
&=(1+t+s)^{\alpha}-(1+s)^{\alpha}\\
&=\alpha\int_s^{t+s}(1+x)^{\alpha-1}\,dx\le\alpha t(1+t+s)^{\alpha-1}
\end{align*}
for $t,s\in\mathbb R^+$, so 
$$\sup_{t\in\mathbb R^+}\frac{|\omega(t+s)-\omega(s)|}{\omega(t)}
\le\frac{\alpha t(1+t+s)^{\alpha-1}}{(1+t)^{\alpha}}
\le\alpha\left(\frac{1+t+s}{1+t}\right)^{\alpha-1}
\le\alpha(1+s)^{\alpha-1}$$ 
for $s\in\mathbb R^+$.

Next, assume that $\alpha<1$. Then the function 
$s\mapsto(1+t+s)^{\alpha}-(1+s)^{\alpha}$ is decreasing on $\mathbb R^+$ for every $t\in\mathbb R^+$, so we deduce that 
$$\omega(t+s)-\omega(s)=(1+t+s)^{\alpha}-(1+s)^{\alpha}
\le(1+t)^{\alpha}-1<\omega(t)$$ 
for $t,s\in\mathbb R^+$. Consequently, 
$$\sup_{t\in\mathbb R^+}\frac{|\omega(t+s)-\omega(s)|}{\omega(t)}\le1$$
for $s\in\mathbb R^+$.

In both cases the result thus follows from Proposition~\ref{pr:cpt3}.
{\nopagebreak\hfill\raggedleft$\Box$\bigskip}

We finish the paper by showing that the condition $\varphi(t)/\omega(t)\to0$ as $t\to\infty$ from Theorem~\ref{th:cpt} cannot in general be relaxed to $\varphi(t)/\omega(t)\to\alpha$ as $t\to\infty$ for some $\alpha\in\mathbb C$.
(For $\varphi\in L^{\infty}(\mathbb R^+)$ we say that $\varphi(t)\to\alpha\in\mathbb C$ as $t\to\infty$ if $\varphi(t)-\alpha\to0$ as $t\to\infty$, that is, if
$\text{ess\,sup}_{t\ge T}\,|\varphi(t)-\alpha|\to0$ as $T\to\infty$.)

\begin{proposition}
\label{pr:notcpt}
Let $\varphi\in L^{\infty}(\mathbb R^+)$ and assume that $\varphi(t)\to\alpha$  as $t\to\infty$ for some $\alpha\neq0$. Then $D_{\varphi}:L^1(\mathbb R^+)\to L^{\infty}(\mathbb R^+)$ and $\overline{D}_{\varphi}:M(\mathbb R^+)\to L^{\infty}(\mathbb R^+)$ are not compact.
\end{proposition}

\noindent{\bf Proof}\quad
Let $\tilde{\varphi}=\varphi-\alpha$, so that $\tilde{\varphi}(t)\to0$ as $t\to\infty$. Let $f_n=1_{[n,n+1]}$ for $n\in\mathbb N$. For $t\in\mathbb R^+$ we have
$$|(D_{\tilde{\varphi}}f_n)(t)|\le\int_n^{n+1}|\tilde{\varphi}(t+s)|\,ds
\le\text{ess\,sup}_{s\ge n}|\tilde{\varphi}(s)|,$$
so we deduce that $D_{\tilde{\varphi}}f_n\to0$ in $L^{\infty}(\mathbb R^+)$ as $n\to\infty$.  
Moreover, 
$$(D_{\alpha}f_n)(t)-\alpha
=\alpha\int_n^{n+1}\frac{s}{t+s}\,ds-\alpha
=-\alpha\int_n^{n+1}\frac{t}{t+s}\,ds$$
for $n\in\mathbb N$ and $t\in\mathbb R^+$, so $\|D_{\alpha}f_n-\alpha\|=|\alpha|$ for $n\in\mathbb N$. Since $D_{\varphi}=D_{\tilde{\varphi}}+D_{\alpha}$, we thus have 
$\|D_{\varphi}f_n-\alpha\|\to|\alpha|$ as $n\to\infty$. 
On the other hand, for $n\in\mathbb N$ and $t\in\mathbb R^+$ we have
$$|(D_{\alpha}f_n)(t)-\alpha|\le|\alpha|\frac{t}{t+n},$$
so it follows from Lebesgue's dominated convergence theorem that $D_{\alpha}f_n\to\alpha$ weak-star in $L^{\infty}(\mathbb R^+)$ and thus $D_{\varphi}f_n\to\alpha$ weak-star in $L^{\infty}(\mathbb R^+)$ as $n\to\infty$. We thus deduce that $(D_{\varphi}f_n)$ has no cluster point as $n\to\infty$. Hence $D_{\varphi}$ and $\overline{D}_{\varphi}$ are not compact.
{\nopagebreak\hfill\raggedleft$\Box$\bigskip}

\bigskip

\noindent
Thomas Vils Pedersen\\
Department of Basic Sciences and Environment\\
University of Copenhagen\\
Thorvaldsensvej 40\\
DK-1871 Frederiksberg C\\
Denmark\\
vils@life.ku.dk


\bigskip
\bigskip

\begin{thebibliography}{10}

\bibitem{Ba-Da:Norms}
{W.\ G.} Bade and {H.\ G.} Dales.
\newblock Norms and ideals in radical convolution algebras.
\newblock {\em J.\ Funct.\ Anal.}, 41:77--109, 1981.

\bibitem{Ba-Da:Cont}
{W.\ G.} Bade and {H.\ G.} Dales.
\newblock Continuity of derivations from radical convolution algebras.
\newblock {\em Studia Math.}, 95:59--91, 1989.

\bibitem{Ch-He:Trans}
Y.~Choi and {M.\ J.} Heath.
\newblock Translation-finite sets and weakly compact derivations from
  {$\ell^1(\Bbb Z_+)$} to its dual.
\newblock {\em Bull.\ Lond.\ Math.\ Soc.}, 42:429--440, 2010.

\bibitem{Ch-He}
Y.~Choi and {M.\ J.} Heath.
\newblock Characterizing derivations from the disk algebra to its dual.
\newblock {\em Proc.\ Amer.\ Math.\ Soc.}, 139:1073--1080, 2011.

\bibitem{Da:Book}
{H.\ G.} Dales.
\newblock {\em Banach algebras and automatic continuity}, volume~24 of {\em
  London Mathematical Society Monographs, New Series}.
\newblock Oxford University Press, Oxford, 2000.

\bibitem{Du-Sc:I}
N.~Dunford and {J.\ T.} Schwartz.
\newblock {\em Linear operators, part {I}}.
\newblock Interscience, New York, 1958.

\bibitem{Gr}
S.~Grabiner.
\newblock Homomorphisms and semigroups in weighted convolution algebras.
\newblock {\em Indiana Univ.\ Math.\ J.}, 37:589--615, 1988.

\bibitem{Gr-wks}
S.~Grabiner.
\newblock Weak{$^*$} properties of weighted convolution algebras {II}.
\newblock {\em Studia Math.}, 198:53--67, 2010.

\bibitem{Gro}
N.~Gr{\o}nb{\ae}k.
\newblock Commutative {B}anach algebras, module derivations, and semigroups.
\newblock {\em J.\ London Math.\ Soc.{(2)}}, 40:137--157, 1989.

\bibitem{Ha}
{P.\ R.} Halmos.
\newblock {\em Measure theory}.
\newblock D.\ Van Nostrand Company, New York, 1950.

\bibitem{Hea}
{M.\ J.} Heath.
\newblock {\em Bounded derivations from {B}anach algebras}.
\newblock PhD thesis, University of Nottingham, 2008.

\bibitem{Hi-Ph}
E.~Hille and {R.\ S.} Phillips.
\newblock {\em Functional analysis and semi-groups}, volume~31 of {\em
  Ame\-rican Mathematical Society Colloquium Publications}.
\newblock Ame\-rican Mathematical Society, Providence, R.\ I., revised edition,
  1957.

\bibitem{Pe:wkprop}
{T.\ V.} Pedersen.
\newblock Weak-star properties of homomorphisms from weighted convolution
  algebras on the half-line.
\newblock {\em J.\ Aust.\ Math.\ Soc.}, 89:75--90, 2010.

\bibitem{Ru}
W.~Rudin.
\newblock {\em Real and complex analysis}.
\newblock McGraw-Hill Book Company, New York, third edition, 1987.

\end{thebibliography}
\end{document}